\documentclass[a4paper,11pt]{article}
\usepackage[utf8]{inputenc}
\usepackage{amsfonts}
\usepackage{amsmath,amsthm}
\usepackage{graphicx}
\usepackage{amssymb}
\usepackage[scale=0.8]{geometry}
\usepackage{tikz}
\usepackage{tikz-cd}
\usepackage{color}
\usepackage{float}
\usepackage{enumerate}
\usepackage{algorithmic}
\usepackage{pifont}
\usepackage[colorlinks,linkcolor=blue,urlcolor=blue,citecolor=blue]{hyperref}
\usepackage{cleveref}
\usepackage{bbold}
\usepackage{mathtools}

\usepackage{dsfont}
\usepackage{ulem}

\newcommand{\calX}{\mathcal{X}}

\newcommand{\RR}{\mathbb{R}}

\newtheorem{theorem}{Theorem}[section]
\newtheorem{remark}{Remark}		
	
\newtheorem{algorithm}[theorem]{Algorithm}	

\newcommand{\cC}{\ensuremath{\mathcal{C}}}

\newcommand{\cP}{\ensuremath{\mathcal{P}}}

\newcommand{\cU}{\ensuremath{\mathcal{U}}}

\newcommand{\cX}{\ensuremath{\mathcal{X}}}

\newcommand{\bN}{\ensuremath{\mathbb{N}}}

\newcommand{\bR}{\ensuremath{\mathbb{R}}}

\newcommand{\xid}{\xi^\dagger}
\newcommand{\xmc}{\bar{x}}

\newcommand{\train}{\ensuremath{\text{train}}}

\DeclareMathOperator*{\argmin}{arg\,min}

\newcommand{\dx}{\ensuremath{\mathrm dx}}

\newcommand{\eps}{\ensuremath{\varepsilon}}

\newcommand{\cond}{\; :\;}

\usepackage{dsfont}

\usepackage{ifthen}
\newlength{\leftstackrelawd}
\newlength{\leftstackrelbwd}
\def\leftstackrel#1#2{\settowidth{\leftstackrelawd}%
{${{}^{#1}}$}\settowidth{\leftstackrelbwd}{$#2$}%
\addtolength{\leftstackrelawd}{-\leftstackrelbwd}%
\leavevmode\ifthenelse{\lengthtest{\leftstackrelawd>0pt}}%
{\kern-.5\leftstackrelawd}{}\mathrel{\mathop{#2}\limits^{#1}}}

\begin{document}

\title{A Dynamical Neural Galerkin Scheme \\for Filtering Problems}

\author{Joubine Aghili
\footnote{IRMA, Université de Strasbourg, CNRS UMR 7501, 7 rue René Descartes, 67084 Strasbourg, France},
Joy Zialesi Atokple\footnote{Kwame Nkrumah University of Science and Technology, Kumasi, Ghana},
Marie Billaud-Friess\footnote{Centrale Nantes, Nantes Université, LMJL UMR CNRS 6629, Nantes, France},
Guillaume Garnier\footnote{Sorbonne Université, Inria, Laboratoire Jacques-Louis Lions (LJLL), 75005 Paris, France},\\
Olga Mula\footnote{Eindhoven University of Technology, Den Dolech 2, P.O. Box 513, 5600 Eindhoven, Netherlands}, 
Norbert Tognon\footnote{INRIA-Paris Project Team ANGE, Sorbonne Université (LJLL), 2 Rue Simone Iff, 75012 Paris, France
}}


\date{\today}
\maketitle

\begin{abstract}
  This paper considers the filtering problem which consists in reconstructing the state of a dynamical system with partial observations coming from sensor measurements, and the knowledge that the dynamics are governed by a physical PDE model with unknown parameters. We present a filtering algorithm where the reconstruction of the dynamics is done with neural network approximations whose weights are dynamically updated using observational data. In addition to the estimate of the state, we also obtain time-dependent parameter estimations of the PDE parameters governing the observed evolution. We illustrate the behavior of the method in a one-dimensional KdV equation involving the transport of solutions with local support. Our numerical investigation reveals the importance of the location and number of the observations. In particular, it suggests to consider dynamical sensor placement.
\end{abstract}


\maketitle

\section{Introduction}
Data assimilation of dynamical systems is a very active research topic, with many challenging open problems especially in cases where strong transport effects are present. This work presents a filtering algorithm based on neural network approximations whose weights dynamically evolve in time.

\subsection{The filtering problem}
We define filtering as the following task. Let $\calX \subseteq \RR^d$ be a spatial domain, and let $V$  a Hilbert space of functions over $\cX$.
In this work, we assume that $V \subset \cC(\cX)$ for simplicity (a more general theoretical framework will be developed in a subsequent contribution). The filtering task is to recover at each time $t\geq0$ the unknown state $u\in \cC^1(\bR_+, V)$ from $m$ observations
$$
z_i(t) = u(t,x_i), \quad 1\leq i \leq m,
$$
where $x_i\in \cX$ are observation points. 

Recovering $u$ only from the observations is a very ill-posed problem so we add an a priori assumption on the dynamics of the system: we assume that $u$ is the solution of an evolution problem of the form
\begin{equation}
\left\{
\begin{array}{rcll}
\partial_t u(t,x,\xi) &=& f\left(t,x,u(t, x, \xi),\xi\right), &(t,x,\xi) \in \bR_+^* \times \calX \times \mathcal{P}, \\
 u(0,x,\xi) &=& u_0(x, \xi) &x \in  \calX, \\
\end{array}
\right.
\label{eq:initialpb}
\end{equation}
with initial condition $u_0: \calX \times \cP \to \RR$, and $f : \bR_+^*  \times \calX \times \bR \times \mathcal{P} \to \RR $ being a nonlinear differential operator. The vector $\xi=(\xi_1, \dots, \xi_p)$ is a vector of scalar parameters ranging in a compact domain $\cP\subset \bR^p$. These parameters can represent physical quantities such as velocity, diffusivity, viscosity,  source terms, or initial conditions. In full generality, $\xi$ could depend on time. However, in the following we will rather think of it as a constant. We assume that problem \eqref{eq:initialpb} is well posed in the sense that, for every $\xi\in \cP$ the problem admits a unique solution $u=u(\xi) \in \cC^1(\bR_+, V)$. In this form, we will regard $u$ as a parameter to solution map $\xi\in \cP\to u(\xi)\in \cC^1(\bR_+, V)$ so that $u(\xi)$ is a function of $t$ and $x$. To simplify notation, we will sometimes write $u(\xi)(t, x)=u(t, x, \xi)$.

In a filtering problem, we assume that there is a parameter $\xi^\dagger$ corresponding to the observed evolution. This implies that we are searching to approximate $u(\xi^\dagger)(t, x)$ for all $t\geq 0$ and $x\in \cX$. However, $\xi^\dagger$ is assumed to be unknown so the dynamics cannot be recovered with traditional forward solvers that compute a solution to problem \eqref{eq:initialpb} with $\xi = \xi^\dagger$.
Instead, our data are the observations $z(t)$ or its first derivatives together with the assumption that the dynamics can be expressed in the form of \cref{eq:initialpb}.

\subsection{Contribution: dynamical reconstruction with neural network approximations}
In this work, we build a filtering algorithm where, for each $t\geq 0$, $u(\xi^\dagger)(t, \cdot)$ is approximated by a feed-forward neural network model. Expressed in more concrete terms: viewing $u(\xi^\dagger)(t, \cdot)$ as a function from $\cX$ to $\bR$, we consider a set $\cU$ of admissible functions from $\cX$ to $\bR$. The functions of $\cU$ are generated by parametrized functions of the form 
\begin{align}
  U: \cX\times \Theta &\to \bR \\
(x, \theta)&\mapsto U(x, \theta)
\end{align}
where $\theta$ is a set of parameters belonging to a space $\Theta$. The set $\cU$ is thus defined as
$$
\cU \coloneqq \{ U(\cdot, \theta) \cond \theta\in \Theta\}.
$$
In this work we focus on the particular case where $\cU$ is generated by feed-forward neural networks with a given architecture. As a simple example, shallow neural networks with inputs in $\cX\subset \bR^d$ and outputs in $\bR$ are defined as follows. For a given $n\in \bN$, let $\theta=(\theta_1,\dots, \theta_n) \in \Theta = (\bR^{d+2})^n \cong \bR^{n(d+2)}$ be our set of parameters with $\theta_i= (c_i,b_i,w_i) \in \bR\times \bR\times \bR^d\cong \bR^{d+2}$.
A shallow neural network with width $n$ is the model $U:\cX \times \Theta \to \bR$ defined as
\begin{equation}
    \label{eq:scalarshallownetwork}
    U(x, \theta)
    :=
    \sum_{i=1}^n c_i \, \phi(x, w_i, b_i),
\end{equation}
where $\phi: \cX \times \bR \times \bR^d \to \bR$ is a smooth activation function, and $w$ and $b$ represent weights and biases respectively. One example for $\phi$ is the ReLU activation unit, defined as
$$
\phi(x, w, b) \coloneqq \max (0, w^Tx+b), \quad\forall x\in \bR^d.
$$
Deeper networks could be defined following similar lines.

The rationale, and originality of our approach lies in the fact that the parameters $\theta$ of our reconstruction will dynamically evolve in time, that is, for every $t\geq 0$, we will find weights $\theta(t)\in \Theta$ such that $U(\cdot, \theta(t))$ approximates $u(\xi^\dagger)(t, \cdot)$, at best given our available information $z(t)$, and the form \eqref{eq:initialpb} of the dynamics.
By allowing the parameter $\theta$ to be time dependent, our strategy allows for  reconstruction as time advances.

One may wonder about the benefit of this dynamical, local in time approach compared to a more global, time-space point of view. In this other point of view, we would typically work with neural networks taking space $x\in \cX$ and time $t\in \bR_+$ as input parameters, and we would search for constant weights $\theta\in \Theta$ in a way that $u(\xi^\dagger)(t, x)$ is approximated well by $U((t, x), \theta)$ for all $(t, x) \in \bR_+\times \cX$ (or a closed time interval $[0, T]$). This time-space approach is less adapted to the filtering problem for several reasons:
\begin{enumerate}
\item We need to wait until the end of the time interval to train the network, and reconstruct in the full time interval.
\item To reconstruct at a single time $t$, we expect to need fewer parameters $\theta(t)$ than in a global time-space approach. This point was confirmed in \cite{Billaud-Friess2017Aug} in the case of forward reduced modeling with linear spaces (instead of nonlinear spaces such as neural networks). It was observed  that the  approximation with dynamically evolving linear spaces was more accurate than a time-space approach involving subspaces with larger dimensions.
\end{enumerate}







\subsection{Connections with earlier works}
Our work connects with a relatively recent research trend aiming at solving inverse problems with neural networks. Probably the field of image analysis was the one of the first to start considering these approximation classes for inversion tasks. In this field, neural networks have produced state of the art results in countless applications and nowadays one could say that they have become the method of choice (see, e.g., \cite{egmont2002image}).

The use of neural networks to address  inverse problems involving PDE physical models is more recent. It is currently attracting a lot of attention, and our work belongs to this trend. Among the relevant contributions that have recently appeared, we may mention \cite{raissi2019physics,bar2019unsupervised, mishra2022estimates, molinaro2023neural}. Our method takes inspiration from the dynamical Neural Galerkin Scheme (NGS) introduced in \cite{Bruna2022Mar}, which uses the same idea for forward PDE simulation. The NGS is actually an internal building block of our strategy as we explain later on.
Compared to most existing works, the main novelty of the NGS introduced \cite{Bruna2022Mar} is that it can be seen as a dynamical (nonlinear) approximation approach. This means that NGS can be derived from some, local in time, Dirac-Frenkel variational principle (see, e.g., \cite{Lubich2008Sep}), related to the set ${\cal U}$ of parametrized functions. In practice, it leads to a system of ODEs for dynamical update of the neural network's weights. Such setting for neural networks particularly extends previous several works done within the context of dynamical low-rank approximation in tensor subsets for forward simulation of PDEs (see, e.g., \cite{Koch2007Apr,Nonnenmacher2008Dec,Lubich2013May,Uschmajew2013Jul}).   

Using dynamical approximations for inverse problems such as filtering is a rather unexplored idea. To the best of our knowledge, only \cite{Lombardi2022, SHNT2023, DY22} and \cite[Chapter 7]{Vidlickova2022} have considered such techniques for linear Kalman filtering. The approximation classes that were used in these works are linear subspaces that are also dynamically updated. The other only existing work using dynamical linear subspaces is \cite{MPV2023}. It comes with additional ideas on structure preservation that are illustrated through the use of symplectic dynamical low rank spaces for Hamiltonian dynamics. \cite{MPV2023} is also the only existing work that considers an algorithm to dynamically update the location of the sensors, and which is based on a rigorous analysis of reconstruction errors. The idea of dynamical sensor placement is important for filtering dynamics involving strong transport effects like the ones we aim to consider here. Compared to these works, the main novelty of this paper is the use of the class of feed-forward neural networks instead of dynamical linear subspaces. Since it is very nonlinear in nature, we expected it to be particularly useful transport dominated dynamics.


\subsection{Organization of the paper}
Section \ref{sec:neural-galerkin} summarizes the NGS of \cite{Bruna2022Mar} for forward PDE simulation. Section \ref{sec:inverse} presents our dynamical filtering algorithm for  filtering. It involves the forward NGS as a building block, and provides both a reconstruction of the state $u$, and a time dependent estimation of the PDE parameters $\xi^\dagger$ associated to the observed trajectory. Finally, Section \ref{sec:numerics} illustrates the behavior of the method in a one-dimensional KdV equation. The study reveals that the method can provide good quality reconstructions provided that the sensors are located inside the support of the solution. Since the location where the solution is supported evolves in time, this means that the locations of observations need to move, and follow the solution dynamics.

\section{The Neural Galerkin Scheme for forward problems}\label{sec:neural-galerkin}

Like for many algorithms for inverse problems and data assimilation, our filtering approach requires to evaluate solutions of \cref{eq:initialpb} for different parameter values $\xi \in \cP$. The task of computing the parameter-to-solution map $\xi \mapsto u(\xi) \in \cC^1(\bR_+, V)$ is usually called the forward problem, as opposed to the inverse problem of estimating $u$ from data observations, which corresponds to the mapping $A:z(t)\mapsto A(z(t))\approx u(\xi^\dagger)(t, \cdot)\in V$. Therefore, as a preparatory step, we explain the scheme that we use for the forward problem, which will serve as a building block for our filtering algorithm. We rely on the NGS that was introduced in \cite{Bruna2022Mar} and which involves a dynamical neural network approximation.

\subsection{The scheme}\label{sec:pdepres}

We carry out the explanation for a fixed, given parameter $\xi\in \cP$ so, to avoid overloading notation, we will omit it in this section. Consequently, we focus on the problem of approximating $u\in \cC^1(\bR_+, V)$ solution to
\begin{equation}
\left\{
\begin{array}{rcll}
\partial_t u(t,x) &=& f\left(t,x,u(t, x)\right), &(t,x) \in \bR_+^* \times \calX, \\
 u(0,x) &=& u_0(x) &x \in  \calX. \\
\end{array}
\right.
\label{eq:initialpb-no-param}
\end{equation}
For every $t\geq 0$, we seek for an approximation of $u(t, \cdot)\in V$ of the form $\mathrm{U}(\cdot, \theta(t))\in \cU\subset V$. 

The NGS from \cite{Bruna2022Mar} consists in deriving an evolution equation for the weights $\theta$ by using the strong residual of the equation \cref{eq:initialpb-no-param}. The main steps are as follows. Assuming that $U$ is differentiable w.r.t.~$\theta$, and $\theta$ differentiable w.r.t.~$t$, by the chain rule we have
\begin{equation}
\partial_t \mathrm{U}(x, \theta(t)) = {\nabla_\theta\mathrm{U}(x, \theta(t))}^T \dot \theta(t),
\end{equation}
where $\nabla_\theta = {(\partial_{\theta_1},\dots,\partial_{\theta_n})}^T \in \RR^n$ stands for the gradient  w.r.t. $\theta$. Thus, for a given $t>0$, the strong residual of the equation at a point $x\in \cX$ is the function $r_t: \Theta\times \dot \Theta\times \cX \to \bR$ defined as
\begin{equation}
r_t(\theta, \eta, x) \coloneqq {\nabla_\theta\mathrm{U}(x, \theta)}^T  \eta - f(t, x, \mathrm{U}(x, \theta)),
\quad \forall (\theta, \eta, x)\in \Theta\times \dot \Theta\times \cX
\end{equation}
where $\dot \Theta$ is the space of time derivatives of $\Theta$. The global residual $R_t : \Theta \times \dot \Theta \to \RR$ is obtained by
computing the squared $L^2_\nu(\cX)$ norm of $x\mapsto r_t(\theta, \eta, x)$, that is,
\begin{align}
R_t(\theta, \eta)
&\coloneqq \frac 1 2\int_\calX |{\nabla_\theta \mathrm{U} (x, \theta)}^T \eta - f(t,x,\mathrm{U}(x, \theta))|^2 \nu(\dx) \\
&= \frac 1 2\Vert {\nabla_\theta \mathrm{U} (\theta,\cdot)}^T \eta - f(t,\cdot,\mathrm{U}(\cdot, \theta)) \Vert^2,
\quad \forall (\theta, \eta)\in \Theta\times \dot \Theta.
\label{eq:global-res}
\end{align}
In the above formula, $\nu\in \cP(\cX)$ is a measure over $\cX$ and $\Vert\cdot\Vert$ is the standard $L^2_\nu(\cX)$ norm associated to the usual scalar product $\langle\cdot,\cdot \rangle$.
In this work, we simply work with the Lebesgue measure $\nu(\dx)=\dx$, and we leave a more in-depth investigation on the optimal choice of $\nu$ for future work.

The NGS seeks to define the curve $\theta: \bR_+\mapsto \Theta$ in such a way that its velocity $\dot\theta(t)$ minimizes the global residual \eqref{eq:global-res} at every $t\in \bR_+^*$. This leads to the variational problem
\begin{equation}
\dot \theta(t) \in \argmin_{\eta \in \dot \Theta} R_t(\theta(t),\eta), \quad \forall t>0.
\label{eq:minJ}
\end{equation}
This local in time optimization problem corresponds to the so-called Dirac-Frenkel variational principle  \cite{Lubich2008Sep}. At the initial time $t=0$, we find $\theta(0) := \theta_0$ by solving the least-squares problem 
\begin{equation}
\theta_0 \in \arg \min_{\theta \in \Theta} 
\frac 1 2 \int_\calX | \mathrm{U} (x, \theta) - u_0(x)|^2 \dx.
\label{eq:min0}
\end{equation}


The necessary optimality conditions of problem \cref{eq:minJ} lead to  the so-called Euler equations characterized by 
\begin{equation}
\nabla_\eta R_t(\theta(t), \dot \theta(t)) = 0, \quad \forall t > 0.
\label{eq:stationnary}
\end{equation}
Since
\[
\nabla_\eta R_t(\theta, \eta)
 = 
\int_\calX \nabla_\theta \mathrm{U} (x, \theta) \left( \nabla_\theta \mathrm{U} {(x, \theta)}^T\eta  - f(t,x,\mathrm{U} (x, \theta)) \right) \dx,
\]
we can write \eqref{eq:stationnary} as
$$
\int_\calX \nabla_\theta \mathrm{U} (x, \theta(t)) \left( \nabla_\theta \mathrm{U} {(x, \theta(t))}^T\dot \theta(t)  - f(t,x,\mathrm{U} (x, \theta(t))) \right) \dx = 0
$$
Note that the above relation can equivalently be written as
\begin{equation}
\label{eq:projTU}
\left<  \nabla_\theta \mathrm{U} {(\cdot, \theta(t))}, r_t(\theta(t), \dot \theta(t), \cdot) \right> =0
\end{equation}
The resulting equation \eqref{eq:projTU} can be formally interpreted as $L^2(\cX)$-projection, for each time $t$, on the tangent space of the ``differential manifold'' ${\cal U}$ at $\mathrm{U}(\cdot, \theta(t))$, which is  the linear space spanned by the partial derivatives with respect to the $\theta_i$ of the non linear map $\mathrm{U}$ (for more detailed discussion on this topic, see, e.g., \cite[Section 2.1]{Berman2023Oct}). 

Defining $M : \Theta \to  \RR^{n \times n}$ and  $F: \bR_+^* \times  \Theta \to \RR^{n}$ as
\begin{equation}
  \begin{aligned}
    M(\theta)  &\coloneqq \int_\calX
    \nabla_\theta \mathrm{U} (x, \theta)\nabla_\theta \mathrm{U} (x, \theta)^T  \dx,\\
    F(t,\theta)  &\coloneqq \int_\calX  \nabla_\theta \mathrm{U} (x, \theta) f(t,x,\mathrm{U} (x, \theta))  \dx,
  \end{aligned}
  \label{def:MF}
\end{equation}
then \cref{eq:stationnary} leads to the following system of Ordinary Differential Equations (ODEs) for the parameters $\theta$
\begin{equation}
\left\{
\begin{array}{rcl}
M(\theta(t)) \dot \theta(t) &=& F(t,\theta(t)) , \quad t>0\\
\theta (0) &=& \theta_0.\\
\end{array}
\right.
\label{eq:paramevol}
\end{equation}

\begin{remark} 
If the initial problem given in \cref{eq:initialpb-no-param} comes with boundary conditions, there exists several possibilities to treat them at the discrete level in \cref{eq:paramevol}.
One can enforce them in a strong way within the model $\mathrm{U}$~\cite[Section 4.4]{Bruna2022Mar}, or weakly by adding conditions in the definition of $\mathrm{U}$, see~\cite[Section 3.2]{Bruna2022Mar}.  
\end{remark}

\subsection{Practical implementation}\label{eq:comp}
To implement the ODE system \eqref{eq:paramevol} in practice, we approximate the integrals over $\calX$ with Monte-Carlo quadrature (MC), and we use classical time integration schemes.
For the sake of simplicity, we present in this section a simple forward Euler time integration, but any other time integration scheme could be considered.
In fact, in our implementation, we have used an explicit Runge-Kutta scheme of order 4.

Let  $0=t_0<t_1<t_2<\cdots<t_k<\dots$ with (possibly non uniform) time steps $\delta t_k\coloneqq t_{k+1}- t_{k} >0$, and suppose we compute at each $t_k$ an approximation  $\theta_k \approx \theta (t_k)$ using a forward Euler scheme.
Thus, starting from $\theta_0$ satisfying \cref{eq:min0}, we are led to the scheme 
\begin{equation}
\tilde M(\theta_k)  \theta_{k+1}  = \tilde M(\theta_{k})   \theta_{k} + \delta t_{k}\tilde F(t_{k+1},\theta_{k}), \quad k\ge 0,\\
\label{eq:EulerScheme}
\end{equation}
where, for every $(t, \theta)\in \bR_+\times \Theta$, $\tilde M(\theta)$  and $\tilde F(t,\theta)$ are the approximations of $M(\theta)$ and $F(t,\theta)$ by MC integration involving $J$ independent uniformly distributed random samples $\{\xmc_j\}_{j=1}^J$ over $\cX$, we have
\begin{equation}
  \begin{aligned}
    \tilde M(\theta)&\coloneqq \frac{1}{J} \sum_{j=1}^J
    \nabla_\theta \mathrm{U} (\xmc_j, \theta)\nabla_\theta \mathrm{U} (\xmc_j, \theta)^T,\\
    \tilde F(t, \theta) &\coloneqq \frac{1}{J} \sum_{j=1}^J \nabla_\theta \mathrm{U} (\xmc_j, \theta) f(t,\xmc_j,\mathrm{U} (\xmc_j, \theta)).
  \end{aligned}
 \label{eq:MFJ}
\end{equation}
In practice, $\tilde M_k\coloneqq \tilde M(\theta_k)$ is a $n\times n$ matrix of the form $\tilde M_k = \frac1J \boldsymbol{V} \boldsymbol{V}^T$. Here, $\boldsymbol{V}$ is a matrix in $\RR^{n \times J }$ whose $j$-th columns corresponds to $ \nabla_\theta \mathrm{U} (\xmc_j, \theta)$.
Similarly, $\tilde F_k\coloneqq \frac1J \boldsymbol{V} \boldsymbol{F}$ is a vector of size $n$ where $\boldsymbol{F}=(F_j)_{1\leq 1\leq J} \in \bR^J$ has entries $F_j=f(t_{k+1},\xmc_j,\mathrm{U} (\xmc_j, \theta_k))$.

In our numerical tests, the matrix $\tilde M_k$ quickly becomes singular due to possible linear dependence of the partial derivatives w.r.t.~$\theta_i$ of $U$. In this work, we address this problem by regularizing $\tilde M_k$ by adding a small perturbation $\varepsilon I$ to $\tilde M_k$ to guarantee invertibility of the matrix. This extra term makes that $\lambda_{\min}(\tilde M_k + \eps I)\geq \eps>0$, which ensures invertibility, and helps to prevent conditionning issues. Alternatively, one can consider the approach described in~\cite{Berman2023Oct}, where the idea is to update only randomly chosen entries of the neural network parameters $\theta$ at each step $k$ of \cref{eq:EulerScheme} to tackle this problem.

\begin{remark}
For some problems (e.g., arising from hyperbolic or transport dominated PDEs) uniformly distributed random samples $\{\xmc_j\}_{j=1}^J$ over $\cX$ may not be efficient to capture dynamical and localized features of the approximation $\mathrm{U}$. In that case, one can rather consider adapted sampling strategies to well capture $\mathrm{U}$ based on the spatial gradient of $\mathrm{U}$ (see \cite{Bruna2022Mar} for details).
\end{remark}

We finish this section by recovering our dependence on a parameter $\xi$ in the forward problem. Note that its presence simply implies that all the above quantities become dependent on it: the neural network parameters $\theta(t)$ will depend on it, and so will $M$, $F$ and their sampled counterparts $\tilde M$ and $\tilde F$ (including $\tilde M_k$ and $\tilde F_k$). In what follows, we record this dependence explicitly in the notation when necessary by adding $\xi$ as a superscript.

\section{A dynamical filtering algorithm}\label{sec:inverse}

In this section, we describe our algorithm for dynamical filtering. 
We recall that the setting here is that we are given observations of the form
$$
z_i(t) = u(\xi^\dagger)(t, x_i), \quad 1\leq i \leq m,
$$
where the $x_i\in \cX$ are observation points, and $\xi^\dagger$ is an unknown parameter. Note that we can leverage $z(t)$ to derive information about the velocity of $u$ at these points since
\begin{equation}
\dot z_i(t) = \lim_{\delta t\to 0} (z_i(t+\delta t)-z_i(t))/\delta t =  \partial_t u(\xi^\dagger)(t, x_i), \quad 1\leq i \leq m.\label{eq:dzi}
\end{equation}
In practice, we are given observations at fixed times such that  $\dot z(t_k)$ can be estimated with finite differences, e.g., $\dot z(t_k) \approx (z(t_{k+1})-z(t_{k}))/\delta t_k$ for $k>0$.

Our filtering algorithm requires both $z(t_k)$ and $\dot z(t_k)$ (or an approximation of it).
Since the parameter $\xi^\dagger\in \cP$ is unknown, our approach is based on a joint parameter-state estimation for each time $t_k$.
This procedure corresponds to algorithm \ref{alg:paramestim}, which proceeds in two steps: a state estimation step and a parameter estimation step, denoted by $\mathbf{(S)}$ and $\mathbf{(P)}$ respectively.
First, we use the parameter $\xi_k$ to compute $\theta_{k}$ by making one step of the forward problem with the NGS from section \ref{sec:neural-galerkin}, namely \cref{eq:EulerScheme}.
Second, we use $\dot z(t_{k+1})$ and the state $\mathrm{U}(\cdot, \theta_{k+1})$ to estimate the parameter $\xi_{k+1}$ in such a way that we minimize the residual of the equation at the observation points.


\begin{algorithm}[Main algorithm]\label{alg:paramestim}
~\begin{itemize}
\item \textbf{Initial time ($k=0, t_0=0$):}
  \begin{enumerate}
  \item[$\mathbf{(S)}$] Given $z(0)$, we find
  \begin{equation}
    \label{eq:min0discrete} 
    \theta_0 \in \arg \min_{\theta \in\Theta} 
    \left\{ \sum_{i=1}^m|z_i(0) - \mathrm{U}(x_i, \theta)|^2 \right\}.
  \end{equation}
  and approximate $u(\xi^\dagger)(0, \cdot) \approx U(\cdot, \theta_0 )$. If the favorable case where the initial condition $u_0$ is fully known, we still need to compute its approximation in $\cU$ so we compute
    \begin{equation}
      \theta_0 \in \arg \min_{\theta \in\Theta} 
      \left\{ \sum_{j=1}^J|u_0(\xmc_j,\xid) - \mathrm{U}(\xmc_j, \theta)|^2 \right\},\label{eq:min0.J}    
    \end{equation}
  \item[$\mathbf{(P)}$] Given $\dot z(0)$ and $\theta_0$, we estimate the parameter $\xi_0$ as
  \[
    \xi_{0} \in \arg\min_{\xi \in \cP}
    \left\{ \displaystyle\sum_{i = 1}^{m} | f(t_1,x_i,\mathrm{U}(x_i, \theta_{0}), \xi) - \dot z_i(0)|^2 \right\}.
  \]
  \end{enumerate}
  
\item \textbf{Inductive step ($t_k\to t_{k+1}$ for $k\geq0$):} Assume we have computed $\mathrm{U}(\cdot, \theta_k)$ and $\xi_k$ for time $t_k$. We now proceed to time $t_{k+1}$. Using the observations $z(t_{k})$ and $\dot z(t_{k})$, we estimate $\theta_{k+1}$ and $\xi_{k+1}$:
\begin{enumerate}
\item[$\mathbf{(S)}$] We update $\theta_{k+1}$ as 
\[
\theta_{k+1}=\theta_k + \delta t_k (\tilde M^{\xi_k}_k)^{-1}\tilde F_k^{\xi_k}. 
\]
We then approximate $u(\xi^\dagger)(t_{k+1}, \cdot) \approx \mathrm{U}(\cdot, \theta_{k+1})$.
\item[$\mathbf{(P)}$] We find $\xi_{k+1}$ as
\begin{equation}
\label{eq:Pstep}
\xi_{k+1} \in \arg\min_{\xi \in \cP}
\left\{ \displaystyle\sum_{i = 1}^{m} | f(t_{k+1},x_i,\mathrm{U}(x_i, \theta_{k+1}), \xi) - \dot z_i(t_{k})|^2 \right\}.
\end{equation}
\end{enumerate}
\end{itemize}
\end{algorithm}

\begin{remark}One could consider more elaborate strategies where the observation locations $x_i$ evolve in time. This feature is important to consider when working with dynamics involving strong transport effects. For these problems, appropriate locations for one time may become completely uninformative for other times due to transport as is shown in \cite{MPV2023}. Even if a systematic strategy on how to place observation points goes beyond the scope of this paper, we explore this feature in our numerical experiments by moving the observations points according to ``pre-scheduled'' dynamics.
\end{remark}

\begin{remark}
  Instead of estimating $\theta_k$ and $\xi_k$ in a two-stage approach, one could alternatively consider a joint estimation of $(\theta_k, \xi_k)$. This point is also left for future work.
  \end{remark}

\section{Numerical experiments}\label{sec:numerics}
  
In what follows, we apply the \cref{alg:paramestim} to treat a filtering problem involving the one-dimensional Korteveg de Vries (KdV) equation \cite{Taha84} presented in~\cite[Section 3.1]{Bruna2022Mar}. For simplicity, we have assumed that $u_0$ was fully known in the experiments so we use \cref{eq:min0.J} in our initialization (instead of \eqref{eq:min0discrete}).
An important aspect which we investigate is the impact of the number and location of the observation points in the reconstruction. We explore the potential benefits of dynamically moving the observation points giving us the observations $z_i(t)$ compared to fixed locations.

\subsection{Problem setting}\label{sec:cont}

\begin{figure}[H]
\centering
\includegraphics[height=140pt]{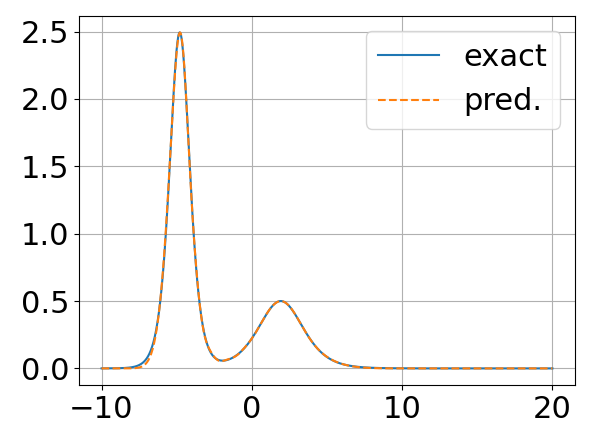}
\includegraphics[height=140pt]{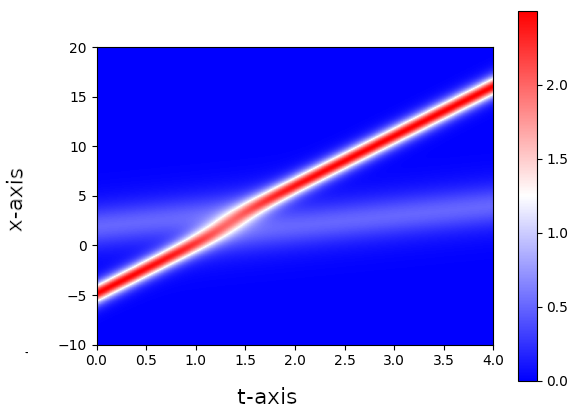}
\caption{\textbf{Left:} Initial condition $u_0$ and its neural network approximation $\mathrm{U}(\cdot, \theta_0)$. \textbf{Right:} Solution to the KdV equation for $\xid=6$. This is the solution that we seek to approximate with the filtering algorithm.}
\label{fig:soln}
\end{figure}

We consider the KdV equation in one space dimension which reads
\begin{equation}
\label{eq:KdV}
\left\{
\begin{array}{ll}
\partial_{t}u(t,x,\xi) &=-\partial_{x}^{3}u(t,x,\xi)-\xi u(t,x,\xi)\partial_{x}u(t,x,\xi),\\
u(0,x,\xi) &= u_0(x),\\
\end{array}
\right.
\end{equation}
where the spatial domain and final time are set to $\mathcal{X}=[-10,20]$ and $T=4$, respectively. 
Here, the initial condition $u_0$ is assumed to be parameter independent. It is defined as in the second case considered in~\cite[pp. 233]{Taha84}. Since its definition is slightly convoluted, we omit giving the full details here and simply plot its shape in \cref{fig:soln} (left plot). The initial condition $u_0$ models the interaction of two solitons with different amplitudes, one located around $x=-5$ and the other around  $x=1$.
Throughout the section, we consider the filtering problem which consists in recovering the solution of \cref{eq:KdV} with fixed $\xi^\dagger=6$.  The map $(x,t)\mapsto u(\xi^\dagger)(t,x)$ is depicted in \cref{fig:soln} (right).

\subsection{Implementation details}\label{sec:disc}

\textbf{Neural network class:} We work with a family $\cU$ of shallow neural networks where $U$ is of the form \eqref{eq:scalarshallownetwork}. Our class consists of a sum of $n=12$ parametrized Gaussian functions corresponding to the choice 
$$
\phi(x,w,b) := \exp (-w^2|x-b|^2).
$$
We therefore have a total of $n(d+2) = 36$ parameters to train (i.e., $\theta\in \bR^{36}$).

\textbf{Forward simulation:}
To solve forward problems, we apply the NGS of section \ref{sec:neural-galerkin} with $\cU$. This leads to the nonlinear system of ODEs given in \cref{eq:paramevol}, and in our case the matrices have size $36\times 36$. We integrate this system with an explicit fourth-order Runge-Kutta (RK4) discretization to produce solutions to the forward problem. Our reference solution from \cref{fig:soln} was obtained in this way.

\textbf{Filtering:} The ODE system \eqref{eq:paramevol} is also solved during the inductive step $\mathbf{(S)}$ of our filtering \cref{alg:paramestim}. We apply the RK4 scheme over $K=1000$ time steps with a fixed time step $\delta t=T/K$. Each step of the scheme requires the evaluations of  the matrices $\tilde{M}$ and the right-hand side $\tilde{F}$ at intermediates time-steps.
They are estimated using $J=1000$ points $\{\xmc_j\}_{j=1}^J$, uniforlmy sampled in $\mathcal{X}$. In the step $\mathbf{(P)}$ of \cref{alg:paramestim}, we replace the continuous interval $\cP$ with $\cP^{\train}\subset \cP$ with $\# \cP^{\train}=101$.

\begin{figure}[H]
  \centering
  \includegraphics[height=140pt]{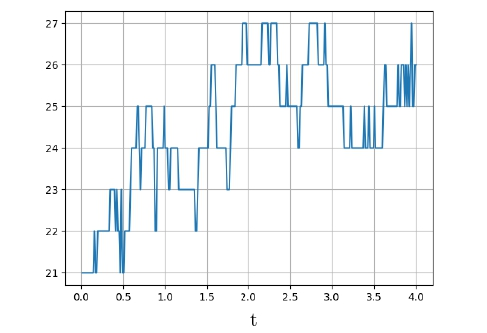}
  \caption{Percentage of eigenvalues of $\tilde M_k$ greater than $10^{-6}$ in absolute value, at each time $t\in[0,4]$.}\label{fig:matrixM}
\end{figure}

To study the conditioning of the matrices $\tilde{M}_k$, we present in 
\cref{fig:matrixM} the percentage of their nonzero eigenvalues for each time step $t_k$ which are larger than $10^{-6}$ in absolute value. We observe that 20--25\% of the eigenvalues are above this treshold.  This indicates that the matrix should not be considered to be invertible, and having full rank. To overcome the rank deficiency, we regularize $\tilde{M}_k$  by adding a regularization term $\eps I$, with $\eps=10^{-3}$.

\textbf{Details about the code:} Our numerical results were obtained using a code written in the \texttt{Python} language. The computation of all derivatives and the minimization routines required for NGS are performed with the \texttt{pytorch} library.

\subsection{Quality of reconstruction}\label{sec:numer}

In this section, we study the reconstruction quality of our filtering algorithm in terms of its ability to approximate $\xi^\dagger$, and $u(\xi^\dagger)$. Remember that $\xi^\dagger$ is constant in time ($\xi^\dagger=6$) but our algorithm computes a time-dependent estimate of it through the $\mathbf{(P)}$ step (see \cref{eq:Pstep}).

To explore the impact of the amount and location of the observation points in the reconstruction, we start by considering in \cref{sec:refsol} a favorable setting with many observation points ($m=100$) uniformly spread in $\cX$. In \cref{sec:impactobs}, we consider a more challenging setting with fewer observations ($m=10$). 

\subsubsection{A favorable setting with many, well-located observations}\label{sec:refsol}

\begin{figure}[H]
  \centering
  \includegraphics[height=140pt]{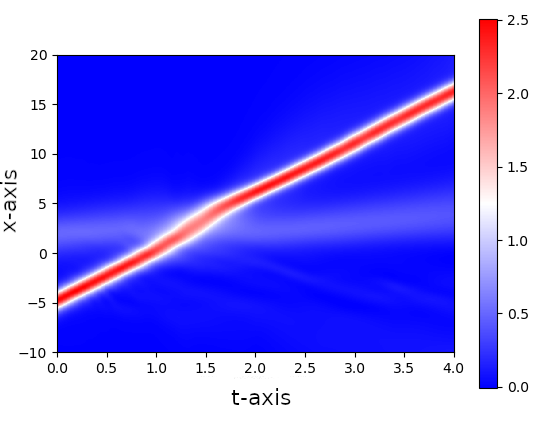} 
  \includegraphics[height=140pt]{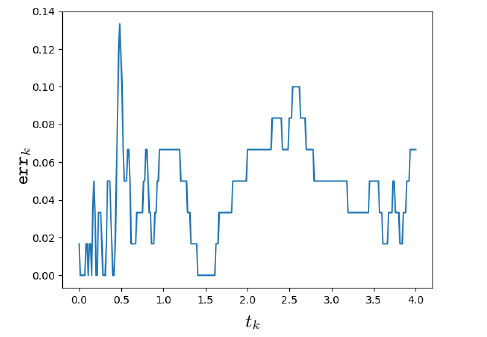} %
  \caption{\textbf{Left:} Reconstructed solution using $m=100$ sensors uniformly spread over $\cX$. \textbf{Right:} Relative error \texttt{err}$_k$ on the approximation of $\xid$.
  }\label{fig:pbinv100}
\end{figure}

\begin{figure}[H]
  \centering
  \includegraphics[width=0.61\paperwidth]{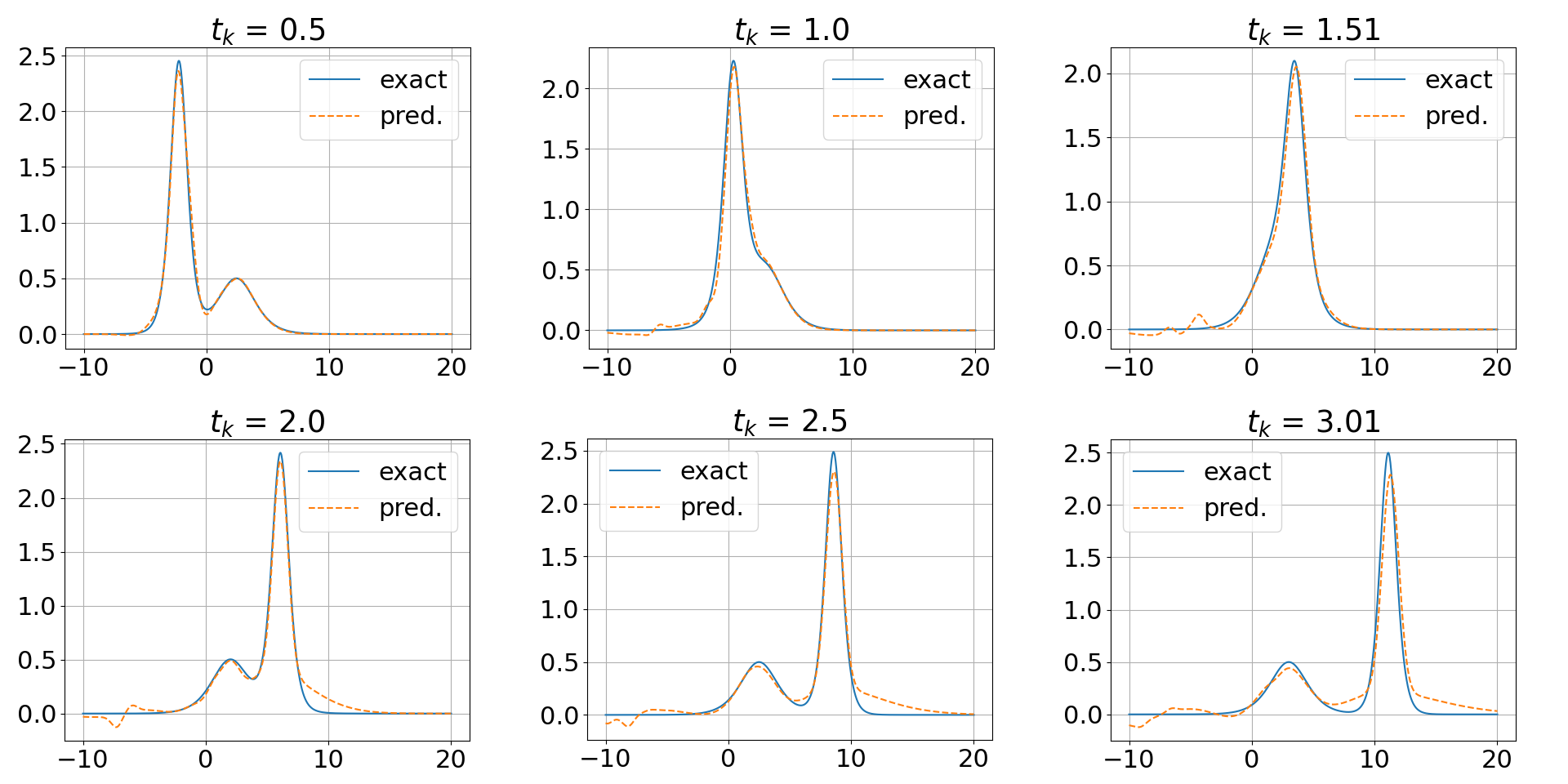}
  \caption{Reconstructed solution using $m=100$ sensors uniformly spread over $\cX$ for different times.}\label{fig:pbinv100t}
\end{figure}

We first choose for $\{x_i\}_{i=1}^m$ a uniform grid of $m=100$ of points in $\cX$.
For our simulations, the quantity $\dot{z}_i$ is estimated using \cref{eq:dzi} with $\delta t=10^{-3}$ from the true solution for $\xid$.
In \cref{fig:pbinv100} we represent the estimated state with respect to time and space.
We also plot the relative error between the true solution and the approximation at each time step $t_k$ as $\mathtt{err}_k:=\frac{|\xid-\xi_k|}{|\xid|}$.
The left figure shows that the estimated state is in very good agreement with the true solution displayed in \cref{fig:soln}. This is confirmed when plotting the reconstructed state for given times $t \in\{0.5,1,1.51,2,2.5,3.01\}$ on \cref{fig:pbinv100t}. 
Similary, the parameter $\xid$ is well estimated, up to an error around 13\%.
These first results demonstrate the ability of our approach to properly estimate both the state $u$ and the parameter $\xid$. The rather good reconstruction results are related to the location and the amount of observation points.
The next section deals with the situation with very few number of observation points.

\subsubsection{A challenging setting with few observations}\label{sec:impactobs}

\begin{figure}[H]
  \centering
  \includegraphics[height=140pt]{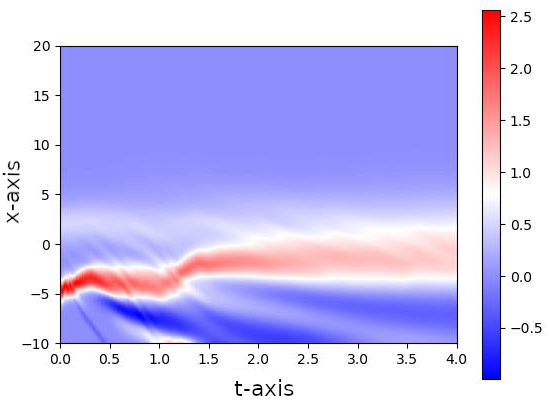}
  \includegraphics[height=140pt]{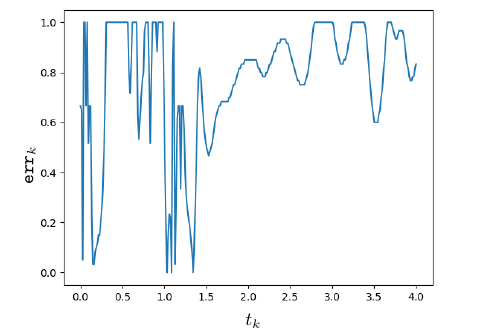}
  \caption{\textbf{Left:} Reconstructed solution using $m=10$ sensors uniformly spread over $\cX$. \textbf{Right:} Relative error \texttt{err}$_k$ on the approximation of $\xid$. }\label{fig:pbinv10}
\end{figure}

\begin{figure}[H]
  \centering
  \includegraphics[width=0.61\paperwidth]{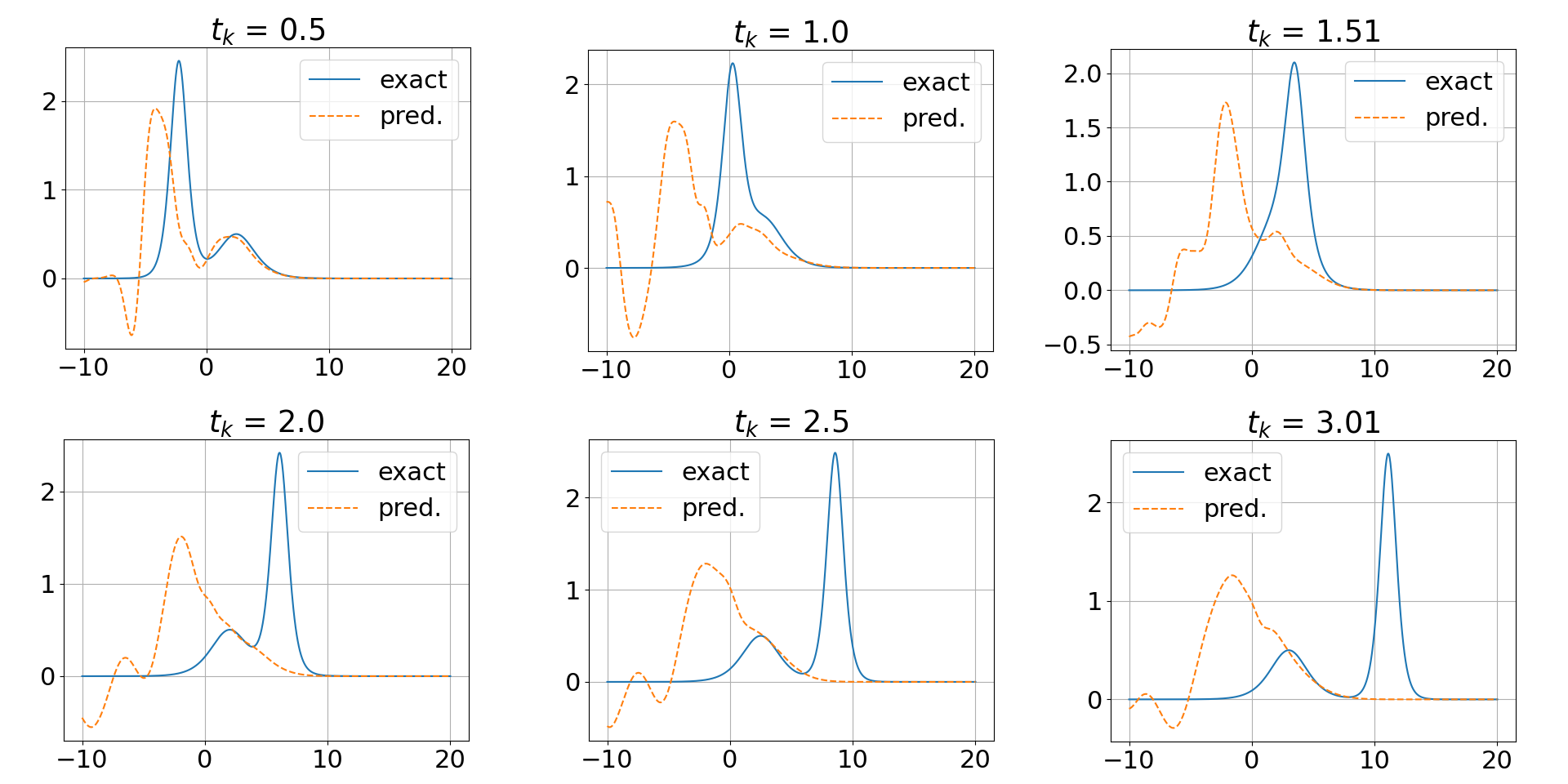}
  \caption{Reconstructed solution  using $m=10$ sensors uniformly spread over $\cX$ for different times.}\label{fig:pbinv10t}
\end{figure}

We start by considering the same setting as above but now we only have $m=10$ observation points uniformly spread over $\cX$. The \cref{fig:pbinv10,fig:pbinv10t} show  the evolution of the reconstructed solution and $\mathtt{err}_k$ of the estimated parameter.
In comparison to \cref{fig:soln}, the dynamics of the solution is clearly not well estimated. This also holds for the parameter estimation for which the relative error can reach 100\%. This behavior is not surprising given that the amount of measurements is rather limited, and their locations only allow to sense the solution very partially at every time.

\begin{figure}[H]
  \centering
  \includegraphics[height=140pt]{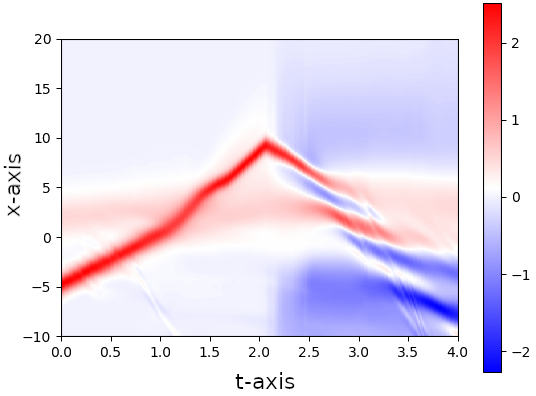}
  \includegraphics[height=140pt]{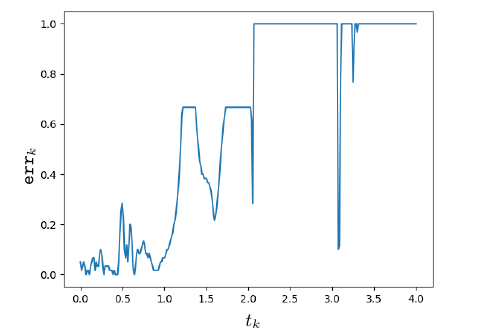}
  \caption{\textbf{Left:} Reconstructed solution using  $m=10$ sensors concentrated around the support of the initial condition. \textbf{Right:} Relative error \texttt{err}$_k$ on the approximation of $\xid$. }\label{fig:pbinv10subset}
\end{figure}

\begin{figure}[H]
  \centering
  \includegraphics[width=0.61\paperwidth]{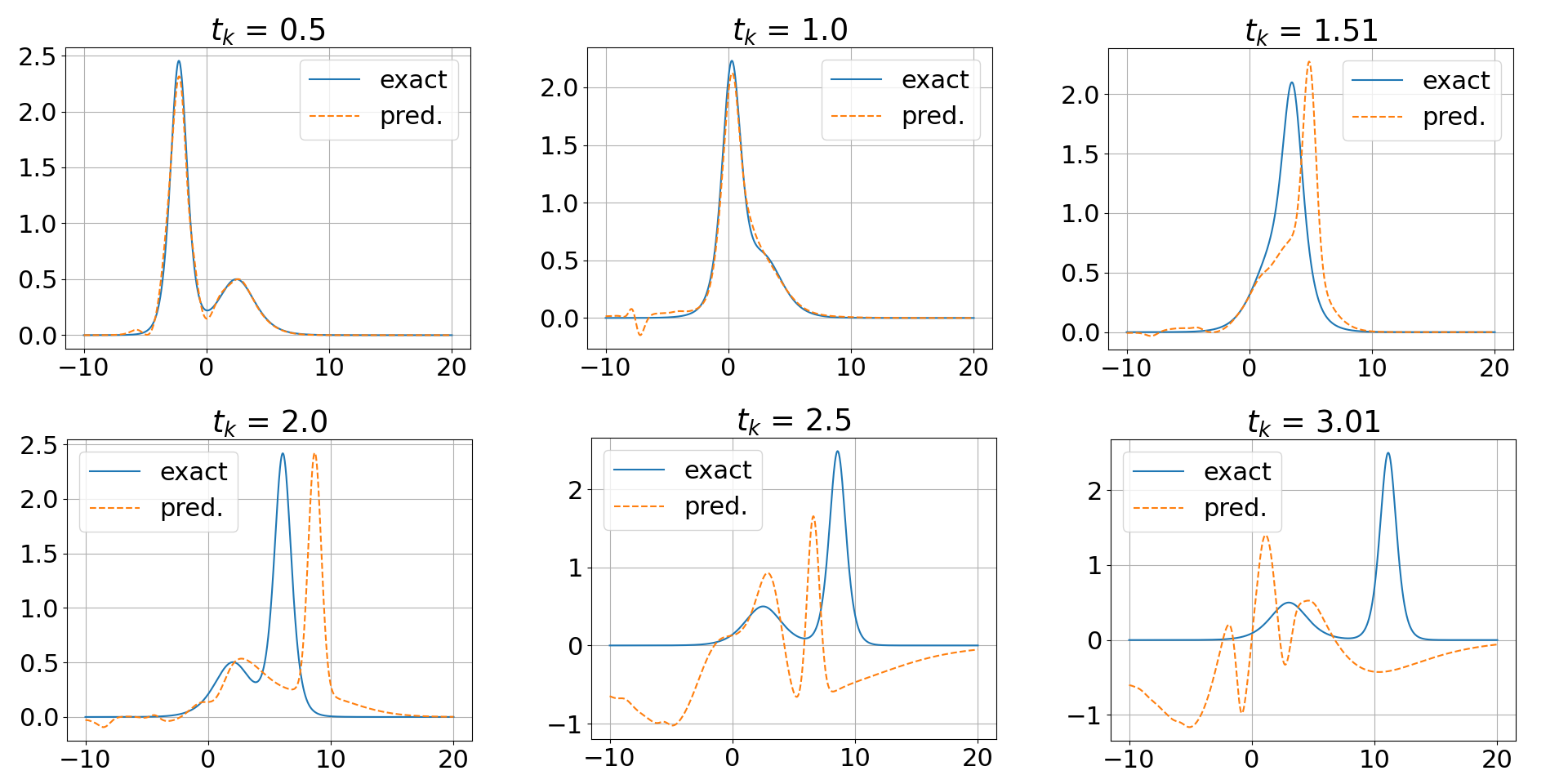}
  \caption{Reconstructed solution using  $m=10$ sensors concentrated around the support of the initial condition  for different times.}\label{fig:pbinv10subsett}
\end{figure}

We next concentrate the $m=10$ sensors around the support of the initial condition $u_0$: we place them uniformly in $[-5, 0]\subset \cX$. One may expect that this location allows for a better reconstruction at early times, and degrades as time advances. The obtained behavior is reported on \cref{fig:pbinv10subset}, and confirms this intuition: for $t\leq 1$, we obtain fairly decent reconstructions with an error up to 30\%. The quality degrades quickly at later times $t>1$. This  is also illustrated on \cref{fig:pbinv10subsett} showing that the solution is poorly reconstructed after $t=1$.

\begin{figure}[H]
  \centering
  \includegraphics[height=140pt]{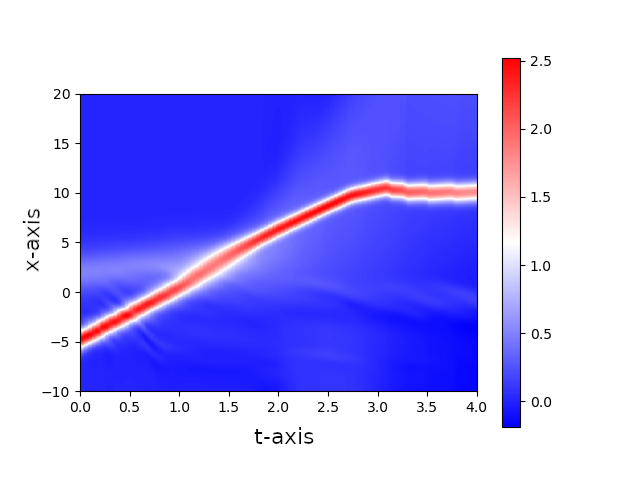}
  \includegraphics[height=140pt]{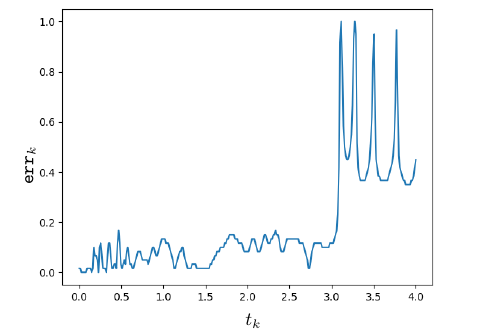}
  \caption{\textbf{Left:} Reconstructed solution  using  $m=10$ sensors that follow the main soliton. \textbf{Right:} Relative error \texttt{err}$_k$ on the approximation of $\xid$.}\label{fig:pbinv.moving}
\end{figure}

\begin{figure}[H]
  \centering
  \includegraphics[width=0.61\paperwidth]{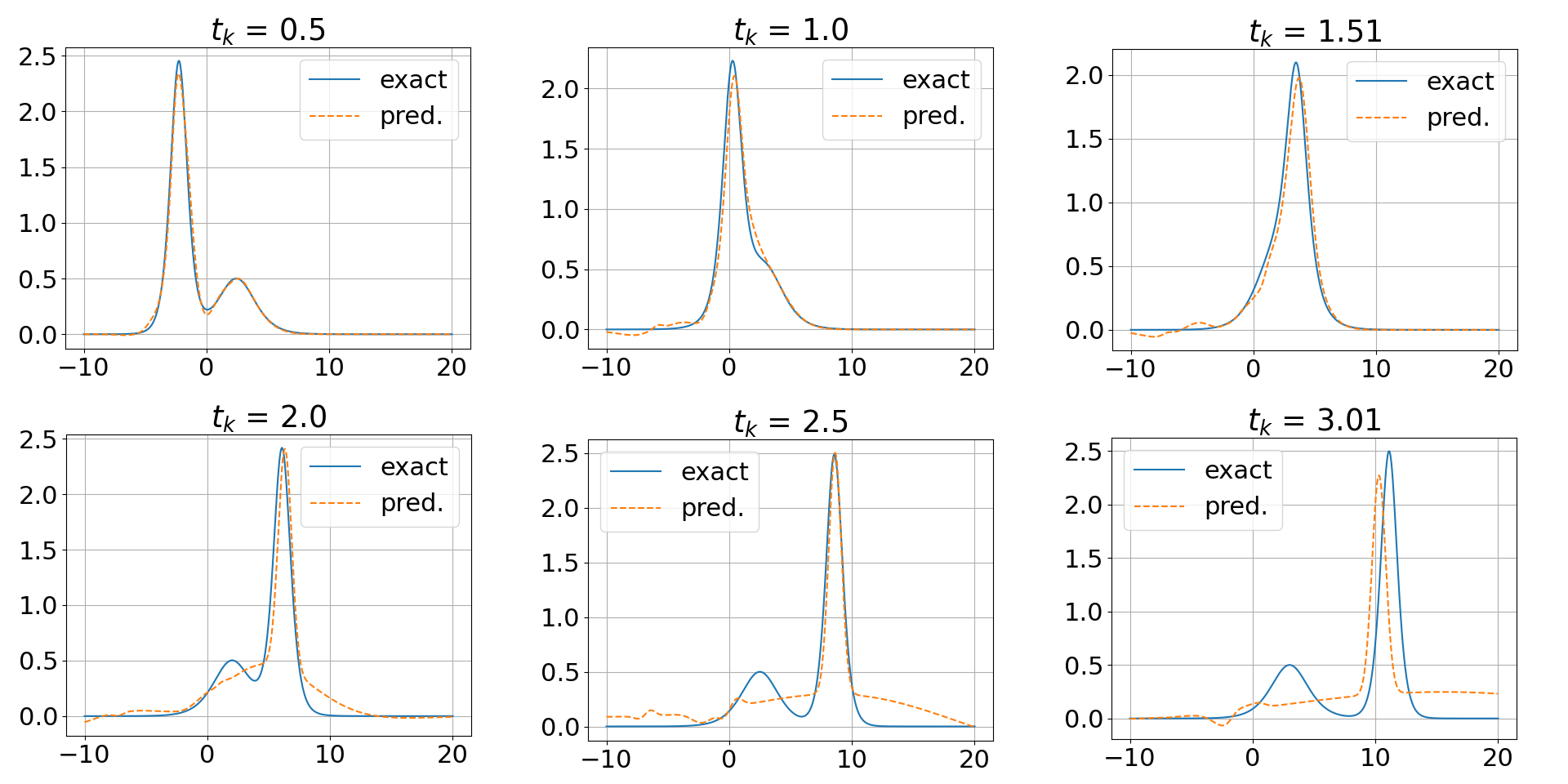}
  \caption{Reconstructed solution  using  $m=10$ sensors that follow the main soliton for different times.
  }\label{fig:pbinv.movingt}
\end{figure}

The above experiments illustrate that, when measurements are scarce, the quality of the estimations highly depends on the localization of the observation points with respect to the support of the solution.
To improve the quality of the estimates, we now consider $m=10$ dynamical and equidistant observation points lying in the time dependent interval $[-5+4t, 0+6t]$. Although this choice is ``hand-made'', it is motivated by the idea of following the location of the main soliton. In full generality, the location should be given by an algorithm that gives the sensors' movement in the course of the evolution in a similar spirit as was done in  \cite{MPV2023} for a filtering algorithm with dynamically updated subspaces. The reconstructed solution over space and time is represented on  \cref{fig:pbinv.moving}.
We observe significant improvements of the estimates in a much longer period of time, i.e. over $[0, 3]$ in comparison to the previous case in \cref{fig:pbinv10,fig:pbinv10subset}.
In the state estimate, we observe that the second soliton (left bump) is badly estimated as shown on \cref{fig:pbinv.movingt}, especially after $t>1.5$. 
However, in comparison with \cref{fig:soln} the main soliton (right bump) is well captured until $t=3.0$.
We can notice a poor parameter estimation after $t>3.0$, as we can see in the relative error plot in \cref{fig:pbinv.moving}.
A better choice for the observation points would have been to track simultaneously both solitons according to their respective characteristics.

\section{Conclusion}\label{sec:conclusion}
We have introduced an approach for filtering which involves dynamically updated neural networks. The method involves the NGS approach from~\cite{Bruna2022Mar} for forward simulation. We have validated the method for a problem where the underlying dynamics of the function to estimate are governed by a 1D parameter dependent KdV equation.
The numerical experiments show that the position of the observation points is of high importance to capture good properties of the solution as well as recovering the latent parameter for the reconstruction with the neural network. When the amount of observations is limited (as is often the case in real applications), our results also indicate that we need to consider methods capable of capturing the dynamics of the function to estimate. There are so far no theoretical foundations for both the forward scheme and the filtering algorithm and this is a crucial point that will be addressed in a forthcoming work.



\begin{thebibliography}{10}

\bibitem{bar2019unsupervised}
Leah Bar and Nir Sochen.
\newblock Unsupervised deep learning algorithm for pde-based forward and
  inverse problems.
\newblock {\em arXiv preprint arXiv:1904.05417}, 2019.

\bibitem{Berman2023Oct}
Jules Berman and Benjamin Peherstorfer.
\newblock {Randomized Sparse Neural Galerkin Schemes for Solving Evolution
  Equations with Deep Networks}.
\newblock {\em arXiv}, October 2023.

\bibitem{Billaud-Friess2017Aug}
Marie Billaud-Friess and Anthony Nouy.
\newblock {Dynamical Model Reduction Method for Solving Parameter-Dependent
  Dynamical Systems}.
\newblock {\em SIAM J. Sci. Comput.}, August 2017.

\bibitem{Bruna2022Mar}
Joan Bruna, Benjamin Peherstorfer, and Eric Vanden-Eijnden.
\newblock {Neural Galerkin Scheme with Active Learning for High-Dimensional
  Evolution Equations}.
\newblock {\em arXiv}, March 2022.

\bibitem{DY22}
G.~Donoghue and M.~Yano.
\newblock A multi-fidelity ensemble {K}alman filter with hyperreduced
  reduced-order models.
\newblock {\em Comput. Methods Appl. Mech. Engrg.}, 398:Paper No. 115282, 18,
  2022.

\bibitem{egmont2002image}
Michael Egmont-Petersen, Dick de~Ridder, and Heinz Handels.
\newblock Image processing with neural networks—a review.
\newblock {\em Pattern recognition}, 35(10):2279--2301, 2002.

\bibitem{Koch2007Apr}
Othmar Koch and Christian Lubich.
\newblock {Dynamical Low{-}Rank Approximation}.
\newblock {\em SIAM J. Matrix Anal. Appl.}, April 2007.

\bibitem{Lombardi2022}
D.~Lombardi.
\newblock State estimation in nonlinear parametric time dependent systems using
  tensor train.
\newblock {\em Internat. J. Numer. Methods Engrg.}, 123(20):4935--4956, 2022.

\bibitem{Lubich2008Sep}
Christian Lubich.
\newblock {From Quantum to Classical Molecular Dynamics: Reduced Models and
  Numerical Analysis}, September 2008.
\newblock [Online; accessed 29. Jan. 2024].

\bibitem{Lubich2013May}
Christian Lubich, Thorsten Rohwedder, Reinhold Schneider, and Bart
  Vandereycken.
\newblock {Dynamical Approximation by Hierarchical Tucker and Tensor-Train
  Tensors}.
\newblock {\em SIAM J. Matrix Anal. Appl.}, May 2013.

\bibitem{mishra2022estimates}
Siddhartha Mishra and Roberto Molinaro.
\newblock Estimates on the generalization error of physics-informed neural
  networks for approximating a class of inverse problems for pdes.
\newblock {\em IMA Journal of Numerical Analysis}, 42(2):981--1022, 2022.

\bibitem{molinaro2023neural}
Roberto Molinaro, Yunan Yang, Bj{\"o}rn Engquist, and Siddhartha Mishra.
\newblock Neural inverse operators for solving pde inverse problems.
\newblock {\em arXiv preprint arXiv:2301.11167}, 2023.

\bibitem{MPV2023}
Olga Mula, Cecilia Pagliantini, and Federico Vismara.
\newblock Dynamical approximation and sensor placement for filtering problems.
\newblock {\em arXiv preprint arXiv:2312.12353}, 2023.

\bibitem{Nonnenmacher2008Dec}
Achim Nonnenmacher and Christian Lubich.
\newblock {Dynamical low-rank approximation: applications and numerical
  experiments}.
\newblock {\em Math. Comput. Simul.}, 79(4):1346--1357, December 2008.

\bibitem{raissi2019physics}
Maziar Raissi, Paris Perdikaris, and George~E Karniadakis.
\newblock Physics-informed neural networks: A deep learning framework for
  solving forward and inverse problems involving nonlinear partial differential
  equations.
\newblock {\em Journal of Computational physics}, 378:686--707, 2019.

\bibitem{SHNT2023}
J.~Schmidt, P.~Hennig, J.~Nick, and F.~Tronarp.
\newblock The rank-reduced {K}alman filter: Approximate dynamical-low-rank
  filtering in high dimensions.
\newblock {\em arXiv preprint arXiv:2306.07774}, 2023.

\bibitem{Taha84}
Thiab~R Taha and Mark~I Ablowitz.
\newblock Analytical and numerical aspects of certain nonlinear evolution
  equations. {III}. {N}umerical, {K}orteweg de {V}ries equation.
\newblock {\em Journal of Computational Physics}, 55(2):231--253, 1984.

\bibitem{Uschmajew2013Jul}
Andr{\ifmmode\acute{e}\else\'{e}\fi} Uschmajew and Bart Vandereycken.
\newblock {The geometry of algorithms using hierarchical tensors}.
\newblock {\em Linear Algebra Appl.}, 439(1):133--166, July 2013.

\bibitem{Vidlickova2022}
E.~Vidli\v{c}kov\'{a}.
\newblock {\em Dynamical low rank approximation for uncertainty quantification
  of time-dependent problems}.
\newblock PhD thesis, EPFL, 2022.

\end{thebibliography}

\end{document}